\newcommand{\al}{\alpha}               \newcommand{\be}{\beta}
\newcommand{\ga}{\gamma}               
\newcommand{\de}{\delta}               
\newcommand{\lb}{\lambda}

        \newcommand{\vphi}{\varphi}

\newcommand{\cal}{\mathcal}
           
\newcommand{\calc}{{\cal C}}           \newcommand{\calf}{{\cal F}}

           \newcommand{\cals}{{\cal S}}

\newcommand{\Fix}{{\rm Fix}}

          \newcommand{\sm}{\setminus}
\newcommand{\impl}{\Rightarrow}      \newcommand{\limpl}{\Longrightarrow}
  
\newcommand{\oo}{\infty}

\newcommand{\uto}{\uparrow}

\newcommand{\bk}{\bigskip}             \newcommand{\sk}{\smallskip}

                \def\R+oo{R_+\cup\{\oo\}}

\def\dtends   {\stackrel {\it d}{\longrightarrow}}

\def\Ctends   {\stackrel {\calc}{\longrightarrow}}

\newcommand{\barr}{\begin{array}}         \newcommand{\earr}{\end{array}}
\newcommand{\bcor}{\begin{corollary}}     \newcommand{\ecor}{\end{corollary}}
\newcommand{\beq}{\begin{equation}}       \newcommand{\eeq}{\end{equation}}
\newcommand{\bex}{\begin{example}}        \newcommand{\eex}{\end{example}}
\newcommand{\bit}{\begin{itemize}}        \newcommand{\eit}{\end{itemize}}
\newcommand{\blemma}{\begin{lemma}}       \newcommand{\elemma}{\end{lemma}}
\newcommand{\bproof}{\begin{proof}}       \newcommand{\eproof}{\end{proof}}
\newcommand{\bprop}{\begin{proposition}}  \newcommand{\eprop}{\end{proposition}}
\newcommand{\brem}{\begin{remark}}        \newcommand{\erem}{\end{remark}}
\newcommand{\btab}{\begin{tabular}}       \newcommand{\etab}{\end{tabular}}
\newcommand{\btheorem}{\begin{theorem}}   \newcommand{\etheorem}{\end{theorem}}

\documentclass[reqno]{amsart}

\newtheorem{theorem}{\bf Theorem}
\newtheorem{corollary}{\bf Corollary}
\newtheorem{example}{\bf Example}
\newtheorem{lemma}{\bf Lemma}
\newtheorem{proposition}{\bf Proposition}
\newtheorem{remark}{\bf Remark}

\begin{document}

\title
[Ran-Reurings Theorems in Ordered Metric Spaces]
{RAN-REURINGS THEOREMS \\
IN ORDERED METRIC SPACES}

\author{Mihai Turinici}
\address{
"A. Myller" Mathematical Seminar;
"A. I. Cuza" University;
700506 Ia\c{s}i, Romania
}
\email{mturi@uaic.ro}

% \date{March 21-26, 2011}

\subjclass[2010]{
47H10 (Primary), 54H25 (Secondary).
}

\keywords{
Ordered metric space, contraction, fixed point, convergence.
}

\begin{abstract}
The Ran-Reurings fixed point theorem
[Proc. Amer. Math. Soc., 132 (2004), 1435-1443]
is but a particular case of Maia's
[Rend. Sem. Mat. Univ. Padova, 40 (1968), 139-143].
A functional version of this last result is then provided, 
in a convergence-metric setting.
\end{abstract}

\maketitle

% Section 1
\section{Introduction}
\setcounter{equation}{0}

Let $X$ be a nonempty set.
Take a metric $d(.,.)$ over it;
as well as a self-map $T:X\to X$.
We say that $x\in X$ is a {\it Picard point} (modulo $(d,T)$) if
{\bf i)} $(T^nx; n\ge 0)$ (=the {\it orbit} of $x$) is $d$-convergent, 
{\bf ii)} $z:=\lim_nT^nx$ is in $\Fix(T)$ (i.e., $z=Tz$).
If this happens for each $x\in X$ and
{\bf iii)} $\Fix(T)$ is a singleton, 
then $T$ is referred to as a {\it Picard operator} (modulo $d$); 
cf. Rus \cite[Ch 2, Sect 2.2]{rus-2001}.
For example, such a property holds whenever 
$d$ is {\it complete} and
$T$ is {\it $d$-contractive}; cf. (b04).
A structural extension of this fact
-- when an {\it order} $(\le)$ on $X$ is being added --
was obtained in 2004 by
Ran and Reurings \cite{ran-reurings-2004}.
For each $x,y\in X$, denote
\bit
\item[(a01)]
$x<> y$ iff either $x\le y$ or $y\le x$
(i.e.: $x$ and $y$ are comparable).
\eit
This relation is reflexive and symmetric; 
but not in general transitive.
Call the self-map $T$, {\it $(d,\le;\al)$-contractive} (for $\al> 0$),  if
\bit
\item[(a02)]
$d(Tx,Ty)\le \al d(x,y)$,\ $\forall x,y\in X$, $x\le y$.
\eit
If this holds for some $\al\in ]0,1[$, 
we say that $T$ is {\it $(d,\le)$-contractive}.

% Theorem 1
\btheorem \label{t1}
Let $d$ be complete and $T$ be $d$-continuous.
In addition, assume that $T$ is $(d,\le)$-contractive and
\bit
\item[(a03)]
$X(T,<>):=\{x\in X; x<> Tx\}$ is nonempty
\item[(a04)]
$T$ is monotone (increasing or decreasing)
\item[(a05)]
for each $x,y\in X$, $\{x,y\}$ has lower and upper bounds.
\eit
Then, $T$ is a Picard operator (modulo $d$).
\etheorem

According to many authors
(cf.
\cite{agarwal-el-gebeily-o-regan-2008},
\cite{ciric-mihet-saadati-2009},
\cite{gwozdz-lukawska-jachymski-2009},
\cite{nieto-rodriguez-lopez-2005},
\cite{o-regan-petrusel-2008}
and the references therein),
this result is credited to be the first extension of the classical 1922
Banach's contraction mapping principle \cite{banach-1922}
to the realm of (partially) ordered metric spaces.
Unfortunately, the assertion is not true:
some early statements of this type 
have been obtained two decades ago by
Turinici \cite{turinici-1986},
in the context of quasi-ordered metric spaces.
(We refer to Section 5 below for details).
\sk

Now, the Ran-Reurings fixed point result 
found some useful applications to 
matrix and differential/integral equations.
So, it cannot be surprising that, soon after, many extensions of
Theorem \ref{t1} were provided;
see the quoted papers for details.
It is therefore natural to discuss
the position of Theorem \ref{t1}
within the classification scheme proposed by
Rhoades \cite{rhoades-1977}.
The conclusion to be derived reads (cf. Section 2):
the Ran-Reurings theorem is but a particular case of the 
1968 fixed point statement in 
Maia \cite[Theorem 1]{maia-1968}.
Further, in Section 3, some extensions are given 
for this last result,
in the context of quasi-ordered convergence almost metric spaces.
Some trivial quasi-order variants of these are then discussed
in Section 4; note that, as a consequence of this, 
one gets the  related contributions in the area due 
Kasahara \cite{kasahara-1976}
and
Jachymski \cite{jachymski-1993},
as well as the order type statement in
O'Regan and Petru\c{s}el \cite{o-regan-petrusel-2008}.
Some other aspects will be delineated elsewhere.

% Section 2
\section{Main result}
\setcounter{equation}{0}

Let $(X,d;\le)$ be an ordered metric space; 
and $T:X\to X$, a self-map of $X$.
Given $x,y\in X$, any subset 
$\{z_1,...,z_k\}$ (for $k\ge 2$) in $X$ with
$z_1=x$, $z_k=y$, and [$z_i<> z_{i+1}$,  $i\in \{1,...,k-1\}$]
will be referred to as a {\it $<>$-chain} 
between $x$ and $y$; the class of all these will 
be denoted as $C(x,y;<>)$.
Let $\sim$ stand for the relation over $X$ attached to $<>$ as 
\bit
\item[(b01)]
$x\sim y$ iff $C(x,y;<>)$ is nonempty. 
\eit
Clearly, $(\sim)$ is reflexive and symmetric; because so is 
$<>$. Moreover, $(\sim)$ is transitive; 
hence, it is an equivalence over $X$.

The following variant of Theorem \ref{t1} is our starting point.

% Theorem 2
\btheorem \label{t2}
Let $d$ be complete and $T$ be $d$-continuous.
In addition, assume that  
$T$ is $(d,\le)$-contractive and
\bit
\item[(b02)]
$T$ is $<>$-increasing [$x<> y$ implies $Tz<> Ty$]
\item[(b03)]
$(\sim)=X\times X$
[$C(x,y; <>)$ is nonempty, for each $x,y\in X$].
\eit
Then, $T$ is a Picard operator (modulo $d$).
\etheorem

This result includes Theorem \ref{t1};
because (a04) $\limpl$ (b02), (a05) $\limpl$ (b03).
[For, given $x,y\in X$,  
there exist, by (a05), some $u,v\in X$ with
$u\le x\le v$, $u\le y\le v$. This yields $x<>u$, $u<> y$; 
wherefrom, $x\sim y$]. In addition, it tells us
that the regularity condition (a03) is superfluous.

The remarkable fact to be noted is that Theorem \ref{t2}
(hence the Ran-Reurings statement as well) 
is deductible from the 
Maia's fixed point statement \cite[Theorem 1]{maia-1968}.
Let $e(.,.)$ be another metric over $X$.
Call $T:X\to X$, 
{\it $(e;\al)$-contractive} (for $\al> 0$) when 
\bit
\item[(b04)]
$e(Tx,Ty)\le \al e(x,y)$, $\forall x,y\in X$; 
\eit
if this holds for some $\al\in ]0,1[$, 
the resulting convention will read as: 
$T$ is {\it $e$-contractive}.
Further, let us say that 
$d$ is {\it  subordinated} to $e$ when 
$d(x,y)\le e(x,y)$, $\forall x,y\in X$.
The announced Maia's result is:

% Theorem 3
\btheorem \label{t3}
Let $d$ be complete and $T$ be $d$-continuous.
In addition, assume that
$T$ is $e$-contractive and $d$ is subordinated to $e$.
Then, $T$ is a Picard operator (modulo $d$).
\etheorem

In particular, when $d=e$, 
Theorem \ref{t3} is just the 
Banach contraction principle \cite{banach-1922}.
However, its potential is much more spectacular; as certified by

% Proposition 1
\bprop \label{p1}
Under these conventions, we have
Theorem \ref{t3} $\limpl$ Theorem \ref{t2};
hence (by the above) 
Maia's fixed point result implies Ran-Reurings'.
\eprop

\bproof
Let $\al\in ]0,1[$ be the number in (a02);
and fix $\lb$ in $]1,1/\al[$. We claim that
\beq \label{201}
\mbox{
$e(x,y):=\sum_{n\ge 0} \lb^n d(T^nx,T^ny)< \oo$,
for all $x,y\in X$.
}
\eeq
In fact, there exists from (b03),
a $(<>)$-chain $\{z_1,...,z_k\}$ (for $k\ge 2$) in $X$ with 
$z_1=x$, $z_k=y$.
By (b02),
$T^nz_i <> T^nz_{i+1}$, $\forall n$, $\forall i\in \{1,...,k-1\}$;
hence, via (a02), 
$d(T^nz_i,T^nz_{i+1})\le \al^nd(z_i,z_{i+1})$, 
for the same ranks $(n,i)$.
But then
$$
d(T^nx,T^ny)\le 
\sum_{i=1}^{k-1}d(T^nz_i,T^nz_{i+1})\le 
\al^n \sum_{i=1}^{k-1}d(z_i,z_{i+1}),\ \ \forall n;
$$
wherefrom (by the choice of $\lb$)
$$
\sum_{n\ge 0} \lb^n d(T^nx,T^ny)\le 
\sum_{n\ge 0} (\lb\al)^n\sum_{i=1}^{k-1}d(z_i,z_{i+1})< \oo;
$$
hence the claim.
The obtained map $e:X\times X\to R_+$ is
{\it reflexive} [$e(x,x)=0$, $\forall x\in X$],
{\it symmetric} [$e(y,y)=e(y,x)$, $\forall x,y\in X$]
and 
{\it triangular} [$e(x,z)\le e(x,y)+e(y,z)$, $\forall x,y,z\in X$].
Moreover, in view of
$$
e(x,y)=d(x,y)+\lb e(Tx,Ty)\ge \lb e(Tx,Ty), \ \forall x,y\in X,
$$
$d$ is subordinated to $e$. 
Note that $e$ is {\it sufficient} in such a case 
[$e(x,y)=0$ $\limpl$ $x=y$];
hence, it is a (standard) metric on $X$.
On the other hand, the same relation tells us that $T$ is 
$(e,\mu)$-contractive for $\mu=1/\lb\in ]\al,1[$; 
hence, (by definition), $e$-contractive.
This, along with the remaining conditions of
Theorem \ref{t2}, shows that
Theorem \ref{t3} applies to these data; wherefrom, all is clear.
\eproof

% Section 3
\section{Extensions of Maia's result}
\setcounter{equation}{0}

From these developments, it follows that
Maia's result \cite[Theorem 1]{maia-1968}
is an outstanding tool in the area;
so, the question of enlarging it is of interest.
A positive answer to this, in a convergence-metric setting, 
will be described below.
\sk

Let $X$ be a nonempty set.
Denote by $\cals(X)$, the class of all sequences $(x_n)$ in $X$.
By a (sequential) {\it convergence structure} on $X$ we mean, as in
Kasahara \cite{kasahara-1976},
any part  $\calc$ of $\cals(X)\times X$ with the properties
\bit
\item[(c01)]
$x_n=x, \forall n\in N  \limpl ((x_n);x)\in \calc$
\item[(c02)]
$((x_n);x)\in \calc \limpl ((y_n);x)\in \calc$, for each
subsequence $(y_n)$ of $(x_n)$.
\eit
In this case, $((x_n);x)\in \calc$ writes $x_n \Ctends x$;
and reads: $x$ is the {\it $\calc$-limit} of $(x_n)$.
The set of all such  $x$ is denoted $\lim_nx_n$; 
when it is nonempty, we say that $(x_n)$ is $\calc$-{\it convergent};
and the class of all these will be denoted $\cals_c(X)$.
Assume that we fixed such an object, with
\bit
\item[(c03)]
$\calc$=separated:\ $\lim_nx_n$ is a singleton, for each $(x_n)$ in $\cals_c(X)$;
\eit
as usually, we shall write $\lim_nx_n=\{z\}$ as $\lim_nx_n=z$.
(Note that, in the 
Fr\'{e}chet terminology \cite{frechet-1906},
this condition is automatically fulfilled,
by the specific way of introducing the ambient convergence;
see, for instance, 
Petru\c{s}el and Rus \cite{petrusel-rus-2006}).
Let $(\le)$ be a {\it quasi-order}
(i.e.: reflexive and transitive relation) over $X$;
and take a self-map $T$ of $X$.
The basic conditions to be imposed are
\bit
\item[(c04)]
$X(T,\le):=\{x\in X; x\le Tx\}$ is nonempty
\item[(c05)]
$T$ is $\le$-increasing ($x\le y$ $\limpl$ $Tx\le Ty$).
\eit
We say that $x\in X(T,\le)$ is a {\it Picard point} (modulo $(\calc,\le,T)$) if
{\bf i)} $(T^nx; n\ge 0)$ is $\calc$-convergent, 
{\bf ii)} $z:=\lim_nT^nx$ is in $\Fix(T)$ and $T^nx\le z$, $\forall n$.
If this happens for each $x\in X(T,\le)$ and
{\bf iii)} $\Fix(T)$ is $(\le)$-{\it singleton} [$z,w\in \Fix(T)$, $z\le w$ $\limpl$ $z=w$],
then $T$ is called a {\it Picard operator} (modulo $(\calc,\le)$). 
Note that, in this case, each $x^*\in \Fix(T)$ fulfills
\beq \label{301}
\forall u\in X(T,\le):\ x^*\le u \limpl u\le x^*;
\eeq
i.e.: $x^*$ is $(\le)$-maximal in $X(T,\le)$. 
In fact, assume that $x^*\le u\in X(T,\le)$. 
By i) and ii), $(T^nu; n\ge 0)$
$\calc$-converges to some $u^*\in \Fix(T)$ with 
$T^nu\le u^*$, $\forall n$; hence, $x^*\le u\le u^*$.
Combining with iii) gives  $x^*=u^*$; 
wherefrom $u\le x^*$.

Concerning the sufficient conditions for such a property,
an early statement of this type was established by 
Turinici \cite{turinici-1986}; cf. Section 5.
Here, we propose a different approach,
founded on ascending orbital concepts
(in short: ao-concepts) and almost metrics.
Some conventions are in order.
Call the sequence $(z_n; n\ge 0)$ in $X$, 
{\it ascending} if $z_i\le z_j$ for $i\le j$;
and {\it $T$-orbital} when $z_n=T^nx$, $n\ge 0$, for some $x\in X$;
the intersection of these concepts is just the precise one.
We say that 
{\bf iv)} $(\le)$ is {\it $(ao,\calc)$-self-closed} 
when the $\calc$-limit of each $\calc$-convergent ao-sequence is an upper bound of it,
{\bf v)} $T$ is {\it $(ao,\calc)$-continuous} if
[$(z_n)$=ao-sequence, $z_n\Ctends z$, $z_n\le z$, $\forall n$] 
implies $Tz_n\Ctends Tz$.
Further, by an {\it almost metric} over $X$ we shall mean any map 
$e:X\times X\to R_+$; supposed to be reflexive triangular and sufficient.
This comes form the fact that such an object 
has all properties of a metric, excepting symmetry.
Call the sequence $(x_n)$, {\it $e$-Cauchy} when 
[$\forall \de> 0$, $\exists n(\de)$: $n(\de)\le n\le m$ 
$\limpl$ $e(x_n,x_m)\le \de$].
We then say that 
{\bf vi)} $(e,\calc)$ is {\it ao-complete}, provided
[(for each ao-sequence) $e$-Cauchy $\limpl$ $\calc$-convergent]. 

Let $\calf(R_+)$ stand for the class of all functions 
$\vphi:R_+\to R_+$. 
Denote by 
$\calf_i(R_+)$, the subclass of all increasing
$\vphi\in \calf(R_+)$;
and by 
$\calf_1(R_+)$, the subclass of all
$\vphi\in \calf(R_+)$ with 
$\vphi(0)=0$ and [$\vphi(t)< t$, $\forall t> 0$].
We shall term $\vphi\in \calf(R_+)$, a {\it comparison function} if
$\vphi\in \calf_i(R_+)\cap\calf_1(R_+)$ 
and [$\vphi^n(t)\to 0$, for all $t> 0$].
[Note that $\vphi\in \calf_1(R_+)$ follows from $\vphi\in \calf_i(R_+)$
and the last property;
cf. Matkowski \cite{matkowski-1977};
but, this is not essential for us].
A basic property of such functions (used in the sequel) is 
\beq \label{302}
(\forall \ga> 0), (\exists \be> 0), (\forall t):\  
0\le t< \ga+\be \limpl \vphi(t)\le \ga.
\eeq 
For completeness, we supply a proof of this,
due to  Jachymski \cite{jachymski-1994}.
Assume that the underlying property fails;
i.e. (for some $\ga> 0$):
\bit
\item[]
$\forall \be> 0$, $\exists t\in [0,\ga+\be[$, such that
$\vphi(t)> \ga$ (hence, $\ga < t< \ga+\be$).
\eit
As $\vphi\in \calf_i(R_+)$, this yields $\vphi(t)> \ga$, $\forall t> \ga$.
By induction, we get (for some $t> \ga$) 
$\vphi^n(t)> \ga$, $\forall n$; 
so (passing to limit as $n\to \oo$) $0\ge \ga$, contradiction.

Denote, for $x,y\in X$:
$H(x,y)=\max\{e(x,Tx),e(y,Ty)\}$, 
$L(x,y)=\frac{1}{2}[e(x,Ty)+e(Tx,y)]$, 
$M(x,y)=\max\{e(x,y),H(x,y),L(x,y)\}$.
Clearly,
\beq \label{303}
M(x,Tx)=\max\{e(x,Tx),e(Tx,T^2x)\},\ \forall x\in X.
\eeq
Call the self-map $T$, 
{\it $(e,M;\le;\vphi)$-contractive}  (for $\vphi\in \calf(R_+)$), if
\bit
\item[(c06)]
$e(Tx,Ty)\le \vphi(M(x,y))$, $\forall x,y\in X$, $x\le y$;
\eit
when this holds for at least one comparison function $\vphi$,
the resulting convention reads: $T$ is {\it extended $(e,M;\le)$-contractive}.

% Theorem 4
\btheorem \label{t4}
Suppose that [in addition to (c04)+(c05)],
$T$ is  extended $(e,M;\le)$-contractive
and $(ao,\calc)$-continuous, 
$(e,\calc)$ is ao-complete,
and 
$(\le)$ is $(ao,\calc)$-self-closed. 
Then, $T$ is a Picard operator (modulo $(\calc,\le)$).
\etheorem

\bproof
Let $x^*,u^*\in \Fix(T)$ be such that $x^*\le u^*$.
By the contractive condition, $e(x^*,u^*)=0$; 
wherefrom, $x^*=u^*$; and so, $\Fix(T)$ is $(\le)$-singleton.
It remains to show that each $x=x_0\in X(T,\le)$ is a Picard point
(modulo $(\calc,\le,T)$). Put $x_n=T^nx$, $n\ge 0$; and 
let $\vphi\in \calf(R_+)$ be the comparison function given 
by the extended $(e,M;\le)$-contractivity of $T$.

{\bf I)}
By the contractive condition and (\ref{303}),
$$
e(x_{n+1},x_{n+2})\le \vphi(M(x_n,x_{n+1}))
=\vphi[\max\{e(x_n,x_{n+1}),e(x_{n+1},x_{n+2})\}], \forall n.
$$
If (for some $n$) the maximum in the right hand side is 
$e(x_{n+1},x_{n+2})$, then (via $\vphi\in \calf_1(R_+)$)
$e(x_{n+1},x_{n+2})=0$; so that (as $e$=sufficient) 
$x_{n+1}\in \Fix(T)$; and we are done. 
Suppose that this alternative fails:
$e(x_{n+1},x_{n+2})\le  \vphi(e(x_n,x_{n+1}))$, for all $n$.
This yields (by an ordinary induction)
$e(x_n,x_{n+1})\le \vphi^n (e(x_0,x_1))$, $\forall n$;
wherefrom $e(x_n,x_{n+1})\to 0$ as $n\to \oo$.

{\bf II)}
We claim that $(x_n; n\ge 0)$ is $e$-Cauchy in $X$.
Denote, for simplicity,
$E(k,n)=e(x_k,x_{k+n})$, $k,n\ge 0$.
Let $\ga> 0$ be arbitrary fixed; and $\be> 0$
be the number appearing in (\ref{302});
without loss, one may assume that $\be< \ga$.
By the preceding step, there exists a rank $m=m(\be)$ such that
\beq \label{304}
\mbox{
$k\ge m$ implies $E(k,1)< \be/2< \be< \ga$.
}
\eeq
The desired property follows from the 
inductive type relation 
\beq \label{305}
\mbox{
$\forall n\ge 0$:\ \ [$E(k,n)< \ga+\be/2$, for each $k\ge m$].  
}
\eeq
The case $n=0$ is trivial; while the case $n=1$ is clear, via (\ref{304}). 
Assume that (\ref{305}) is true, 
for all $n\in \{1,...,p\}$ (where $p\ge 1$); 
we want to establish that it holds as well for $n=p+1$.
So, let $k\ge m$ be arbitrary fixed.
By the induction hypothesis and (\ref{304}), 
$e(x_k,x_{k+p})=E(k,p)< \ga+\be/2$ 
and 
$H(x_k, x_{k+p})=\max\{E(k,1), E(k+p,1)\}< \be/2$. 
Moreover, the same premises give (by the triangular property)
$$ \barr{l}
L(x_k,x_{k+p})=(1/2)[E(k,p+1)+E(k+1,p-1)] \le \\
(1/2)[E(k,p)+E(k+p,1)+E(k+1,p-1)]< \ga+ \be;
\earr
$$
wherefrom  $M(x_k,x_{k+p})< \ga+\be$; 
so, by the contractive condition and (\ref{302}),
$$
E(k+1,p)=e(x_{k+1},x_{k+p+1})=e(Tx_k,Tx_{k+p})\le 
\vphi(M(x_k,x_{k+p}))\le \ga;
$$
which "improves" the previous evaluation (\ref{305}) of our quantity.
This, along with (\ref{304}) and the triangular property, gives
$E(k,p+1)=e(x_k,x_{k+p+1})< \ga+\be/2$.

{\bf III)}
As $(e,\calc)$ is ao-complete,
(\ref{305}) tells us that $x_n\Ctends x^*$ for some $x^*\in X$.
Moreover, as $(\le)$ is $(ao,\calc)$-self-closed, 
we have $x_n\le x^*$, $\forall n$;
hence, in particular, $x\le x^*$.
Combining with the  $(ao,\calc)$-continuity of $T$,
yields
$x_{n+1}=Tx_n\Ctends Tx^*$; 
wherefrom (as $\calc$ is separated),
$x^*\in \Fix(T)$.
\eproof

Now letting $d$ be a metric on $X$,
the associated convergence $\calc:=(\dtends)$ is separated;
moreover, the  ao-complete property of 
$(e,\calc)$ is holding whenever $d$ is complete and
subordinated to $e$.
Clearly, this last property is trivially assured if $d=e$;
when Theorem \ref{t4} is comparable with 
the main result in 
Agarwal, El-Gebeily and  O'Regan \cite{agarwal-el-gebeily-o-regan-2008}.
In fact, a little modification of the  working hypotheses
allows us getting the whole conclusion of the quoted statement; 
we do not give details.

% Section 4
\section{Particular aspects}
\setcounter{equation}{0}

Let $X$ be a nonempty set; 
and $T:X\to X$ be a self-map of $X$.
Further, take a separated (sequential)
convergence structure  $\calc$ on $X$.
\sk

{\bf (A)}
Let  $e(.,.)$ be an almost metric over $X$. 
A basic particular case of the previous developments corresponds to 
$(\le)=X\times X$ (=the {\it trivial} quasi-order on $X$).
Then, (c04)+(c05) are holding; 
and the resulting Picard concept becomes a 
Picard property (modulo $\calc$) of $T$, which writes:
{\bf i)}
$\Fix(T)$ is a singleton, $\{x^*\}$,
{\bf ii)} 
$T^n x\Ctends x^*$, for each $x\in X$.
Moreover, 
the $(ao,\calc)$-self-closedness of $(\le)$ is fulfilled;
and the remaining ao-concepts become 
orbital concepts (in short: o-concepts).
Precisely, call $T$, {\it $(o,\calc)$-continuous} if
[$(z_n)$=o-sequence, $z_n\Ctends z$ imply $Tz_n\Ctends Tz$];
likewise, call $(e,\calc)$, {\it o-complete} if
[(for each o-sequence) $e$-Cauchy $\limpl$ $\calc$-convergent]. 
Finally, concerning (c06), 
let us say that $T$ is {\it $(e;M;\vphi)$-contractive} 
(where $\vphi\in \calf(R_+)$), provided
\bit
\item[(d01)]
$e(Tx,Ty)\le \vphi(M(x,y))$, $\forall x,y\in X$;
\eit
if this holds for at least one comparison function $\vphi$,
the resulting convention reads: $T$ is {\it extended $(e,M)$-contractive}.
Putting these together, one gets
the following version of Theorem \ref{t4}:

% Corollary 1
\bcor \label{c1}
Suppose that  
$T$ is extended $e$-contractive
and $(o,\calc)$-continuous, 
and 
$(e,\calc)$ is o-complete.
Then, $T$ is a Picard operator (modulo $\calc$).
\ecor

The obtained statement includes
Kasahara's fixed point principle \cite{kasahara-1976},
when $e$ is a metric on $X$.
On the other hand, if $d$ is a metric on $X$
and $\calc:=(\dtends)$, the  o-complete property of 
$(e,\calc)$ is assured when $d$ is complete and
subordinated to $e$.
This, under a linear choice of the comparison function
($\vphi(t)=\al t$, $t\in R_+$, for $0< \al< 1$), 
tells us that Corollary \ref{c1} includes Theorem \ref{t3}.
Finally, when $d=e$, Corollary \ref{c1} reduces
to Jachymski's result  \cite{jachymski-1993}.
\sk

{\bf (B)}
An interesting version of Corollary \ref{c1}
was provided in the 2008 paper by
O'Regan and Petru\c{s}el \cite[Theorem 3.3]{o-regan-petrusel-2008}.
Let $(X,T,\calc)$ be endowed with their precise general meaning;
and $d(.,.)$ be a (standard) metric on $X$.
As before, we are interested to give sufficient conditions 
under which $T$ be a Picard operator (modulo $\calc$).
Take an {\it order} $(\le)$ on $X$; 
and put $X_{(\le)}=(\le)\cup (\ge)$, 
where $(\ge)$ stands for the {\it dual} order.
This subset is just the graph of the
relation $<>$ over $X$ introduced as in (a01);
so, it may be identified with the underlying  relation.
As a consequence, $X_{(\le)}$ is 
{\it reflexive} [$(x,x)\in X_{(\le)}$, for each $x\in X$]
and
{\it symmetric} [$(x,y)\in X_{(\le)}$ iff $(y,x)\in X_{(\le)}$];
but not in general transitive, as simple examples show.
Further, let us say that $(d,\calc)$ is {\it o-complete} if
[(for each sequence) $d$-Cauchy $\limpl$ $\calc$-convergent]. 
Finally, call $T$, 
{\it $(d,\le;\vphi)$-contractive}  (for $\vphi\in \calf(R_+)$), if
\bit
\item[(d02)]
$d(Tx,Ty)\le \vphi(d(x,y))$, $\forall x,y\in X$, $x\le y$;
\eit
when this holds for at least one comparison function $\vphi$,
the resulting convention reads: $T$ is {\it $(d,\le)$-contractive}.

% Corollary 2
\bcor \label{c2}
Assume that 
(a03)+(b02) hold,
$T$ is $(d,\le)$-contractive and $(o,\calc)$-continuous, 
$(d,\calc)$ is complete, and
\bit
\item[(d03)]
$(x,y), (y,z)\in X_{(\le)}$ $\impl$ $(x,z)\in X_{(\le)}$
(i.e.: $X_{(\le)}$ is transitive)
\item[(d04)]
$(x,y)\notin X_{(\le)}$ $\impl$ $\exists$ 
$c=c(x,y)\in X$:\  $(x,c),(y,c)\in X_{(\le)}$.
\eit
Then, $T$ is a Picard map (modulo $\calc$).
\ecor

\bproof
We claim that Corollary \ref{c1} is applicable to 
such data. This will follow from 
$$
\mbox{
$X\times X=X_{(\le)}$\ (i.e.: the ambient order $(\le)$ is linear).
}
$$
In fact, let $x,y\in X$ be arbitrary fixed.
If $(x,y)\in X_{(\le)}$, we are done; so, assume that
$(x,y)\notin  X_{(\le)}$. By (d04), there exists $c=c(x,y)\in X$
such that $(x,c)\in X_{(\le)}$, $(y,c)\in X_{(\le)}$.
This, along with the symmetry of $X_{(\le)}$, gives 
$(c,y)\in X_{(\le)}$; hence, by (d03), 
$(x,y)\in X_{(\le)}$.
As a consequence, the $(d,\le;\vphi)$-contractive property for $T$
is to be written in the "amorphous" form of (d01)
[with $d(.,.)$ in place of $e(.,.)$ and $M(.,.)$]; 
wherefrom, all is clear.
Note that, from such a perspective, 
conditions (a03)+(b02) are superfluous.
\eproof

% Section 5
\section{Old approach (1986)}
\setcounter{equation}{0}

In the following, a summary of the 1986 results in 
Turinici \cite{turinici-1986} 
is being sketched, for completeness reasons.
\sk 

{\bf (A)}
Let $(X,d)$ be a complete metric space and $T$ be a self-map of $X$. 
Assume that for each $x\in X$ there exists a 
$n(x)\in N_0:=N\sm \{0\}$ such that 
$T^{n(x)}$ is (metrically) contractive at $x$; 
then, we may ask of under which additional conditions 
is $T$ endowed with a Picard property (cf. Section 1).
A first answer to this question was given, in the continuous case, by 
Sehgal \cite{sehgal-1969}
through a specific iterative procedure;  
a reformulation of it for discontinuous maps was  performed in
Guseman's paper  \cite{guseman-1970}.  
During the last decade, some technical extensions 
-- involving the contractive condition --
of these results were obtained by 
Ciric  \cite{ciric-1972}, 
Khazanchi \cite{khazanchi-1974}, 
Iseki \cite{iseki-1974},
Rhoades	\cite{rhoades-1977}
and 
Singh \cite{singh-1979}. 
The most general statement of this kind, obtained by 
Matkowski  \cite{matkowski-1977}, 
reads as follows.
For each $m\in N_0$, 
let $\calf(R_+^m)$ stand for 
the class of all functions $f:R_+^m\to R_+$;
and  $\calf_i(R_+^m)$ the subclass of all $f\in \calf(R_+^m)$,
increasing in each variable.
The iterative contraction property below is considered:
\bit
\item[(e01)]
$\exists f\in \calf_i(R_+^5)$ such that:
$\forall x\in X$, $\exists n(x)\in N_0$ with
$d(T^{n(x)}x,T^{n(x)}y) \le$ \\
$f(d(x,T^{n(x)}x),d(x,y),d(x,T^{n(x)}y),d(T^{n(x)}x,y),d(T^{n(x)}y,y))$,\ $\forall y\in X$.
\eit
Given $f\in \calf_i(R_+^5)$ like before, denote
$g(t)=f(t,t,t,2t,2t)$, $t\ge 0$;
clearly, it is an element of $\calf_i(R_+)$.
We shall say that $f$ is {\it normal} provided
\bit
\item[(e02)]
$g\in \calf_1(R_+)$ and [$t-g(t)\to \oo$  as   $t\to \oo$]
\item[(e03)]
$\lim_ng^n(t)= 0$, for each $t> 0$. 
\eit
(As already remarked, (e03) implies the
first part of (e02), under the properties of $g$;
we do not give details).

% Theorem 5. 
\begin{theorem} \label{t5}
Suppose that there exists a normal function
$f\in \calf_i(R_+^5)$ in such a way that (e01) holds.
Then, $T$ is a Picard operator (modulo $d$).
\end{theorem}

A direct examination of the above conditions shows that, by virtue of 
$$
d(T^{n(x)}x,y)\le d(x,T^{n(x)}x)+d(x,y),\ x,y \in X,
$$
$$
d(T^{n(x)}y,y)\le d(x,T^{n(x)}y)+d(x,y),\ x,y \in X,
$$
a slight extension of Theorem \ref{t5} might be reached 
if one replaces (e01) by 
\bit
\item[(e04)]
$d(T^{n(x)}x,T^{n(x)}y) \le F(d(x,T^{n(x)}x),d(x,y),d(x,T^{n(x)}y))$,\ $y\in X$,
\eit
where $F:R_+^3\to R_+$ is defined as
\bit
\item[]
$F(\xi,\eta,\zeta)=f(\xi,\eta,\zeta,\xi+\eta,\zeta+\eta)$,\  \
 $\xi,\eta,\zeta\in R_+$.    
\eit
A natural question to be solved is that of determining 
what happens when the right-hand side of (e04) depends on 
the (abstract) variable $x\in X$ and the (real) variables  
$((d(x,T^ix); 1\le i\le n(x))$, $(d(x,T^jy); 0\le j\le n(x))$; 
or, in other words, when the function $F=F(x)$ acts from $R_+^{2n(x)+1}$ to $R_+$.
At the same time, observe that, from a "relational" viewpoint, 
the result we just recorded may be deemed as being expressed modulo 
the {\it trivial} quasi-ordering on $X$; so that, a formulation 
of it in terms of  genuine quasi-orderings would be of interest. 
It is precisely our main aim to get a generalization -- 
under the above lines -- 
of Theorem \ref{t5}. 
\sk

{\bf (B)}
Let $(X,d)$ be a metric space and $\le$ be a quasi-ordering 
(i.e.: reflexive and transitive relation) over $X$. A sequence $(x_n; n\in N)$ 
in $X$ will be said to be {\it increasing} when $x_i\le x_j$ for $i\le j$.
Take the self-map $T$ of $X$ according to
\bit
\item[(e05)]
$Y:=\{x \in X; x\le Tx\}$ is not empty
\item[(e06)]
$T$ is increasing ($x\le y$ implies $Tx\le Ty$).
\eit
In addition, the specific condition will be accepted:
\bit
\item[(e07)]
for each $x$ in $Y$ there exist $n(x)\in N_0$,
$f(x)\in \calf_i(R_+^{2n(x)+1})$, with \\
$d(T^{n(x)}x,T^{n(x)}y)\le$ 
$f(x)(d(x,Tx),...,d(x,T^{n(x)}x);d(x,y),...,d(x,T^{n(x)}y))$,\\
for all $y\in Y$ with $x\le  y$.
\eit
For the arbitrary fixed $x\in Y$,  let $g(x)$ indicate 
the element of $\calf_i(R_+)$, given as 
$g(x)(t)=f(x)(t,...,t;t,...,t)$, $t\ge 0$. 
We shall say that the family (of (e07)) $((n(x),f(x)); x\in Y)$
is {\it iterative $T$-normal} provided, for each $x_0\in Y$,
\bit
\item[(e08)]
$g(x_0)\in \calf_1(R_+)$  and $t-g(x_0)(t)\to \oo$  as $t\to  \oo$,
\item[(e09)]
$\lim_kg(x_k)\circ...\circ g(x_0)(t)=0$, $t> 0$, where 
[$n_0=n(x_0)$, $x_1=T^{n_0}x_0$] and, inductively, 
[$n_i=n(x_i)$, $x_{i+1}=T^{n_i}x_i$], $i\ge 1$.
\eit
The following auxiliary fact will be useful.

% Proposition 2.     
\begin{proposition} \label{p2}
Let (e05)-(e07) hold; and the family 
$((n(x),f(x)); x\in Y)$ [attached to (e07)] be iterative $T$-normal.
Then, the following conclusions hold

{\bf i)}  
for each $x \in Y$, $(T^mx; m\in N)$ is increasing Cauchy (in $X$) 

{\bf ii)} 
$d(T^mx,T^my)\to 0$ as $m\to \oo$, for all $y\in Y$, $x\le y$.
\end{proposition}

\begin{proof}
Let $x\in Y$ be given. We firstly claim that
\begin{equation} \label{501}
\mbox{
$d(x,T^mx)\le t$, $m \in N$, for some  $t=t(x)> 0$.
}
\end{equation}
Indeed, it follows by (e08) that, given $\al> 0$, there exists $\be=\be(\al,x)> \al$  with
\begin{equation} \label{502}
\mbox{
$t\le \al+ g(x)(t)$  implies $t\le \be$.
}
\end{equation}
Put $\al=\max\{d(x,Tx),...,d(x,T^{n(x)}x)\}$. 
We claim that (\ref{501}) holds with $t=\be$. 
In fact, suppose that the considered assertion would be false; and let $m$   
denote the infimum of those ranks for which the reverse of (\ref{501}) takes place.
Clearly, $m>n(x)$, $d(x,T^kx)\le \be$, $k\in \{1,...,m-1\}$, and $d(x,T^mx)>\be$; 
so that, by (e07), 
$$  \barr{l}
d(x,T^mx)\le d(x,T^{n(x)}x)+d(T^{n(x)}x,T^mx)\le \\ 
\al + f(x)(d(x,Tx),...,d(x,T^{n(x)}x); d(x,T^{m-n(x)}x),...,d(x,T^mx))\le \\
\al + f(x)(\al,...,\al;\be,...,\be,d(x,T^mx))\le  \al + g(x)(d(x,T^mx));
\earr
$$
contradicting (\ref{502}) and proving our assertion. 
In this case, letting $x= x_0$ in $Y$, put $n_0=n(x_0)$, $m_0=n_0$, 
$x_1=T^{n_0}x_0=T^{m_0}x_0$
and, inductively, 
$$
n_i=n(x_i), m_i= n_0+...+n_i, 
x_{i+1}= T^{n_i}x_i= T^{m_i}x_0,\ i\ge 1. 
$$
By (\ref{501}), $d(x_0,T^mx_0)\le t_0$, $m\in N$, 
for some $t_0> 0$; so combining with (e07):
$$  \barr{l}
d(x_1,T^mx_1)=d(T^{n_0}x_0, T^{n_0}T^mx_0) \le 
f(x_0)(d(x_0,Tx_0),...,d(x_0,T^{n_0}x_0); \\
d(x_0,T^mx_0),...,d(x_0,T^{n_0+m}x_0))\le 
g(x_0)(t_0), m\in N;
\earr
$$
or equivalently, 
$d(T^{m_0}x_0,T^mx_0)\le g(x_0)(t_0)$, $m\ge m_0$.
Again via (e07),  
$$ \barr{l}
d(x_2,T^mx_2)=d(T^{n_1}x_1,T^{n_1}T^mx_1)\le 
f(x_1)(d(x_1,Tx_1),...,d(x_1,T^{n_1}x_1); \\
d(x_1,T^mx_1),...,d(x_1,T^{n_1+m}x_1))\le 
g(x_1)\circ g(x_0)(t_0), m\in N;
\earr
$$
or equivalently:
$d(T^{m_1}x_0,T^mx_0)\le g(x_1)\circ g(x_0)(t_0)$, $m\ge m_1$; and so on. 
By a finite induction procedure one gets
$d(x_{k+1},T^mx_{k+1})\le g(x_k) \circ...\circ g(x_0)(t_0),   m,k\in N$;
or equivalently (for each $k\in N$)
$$
d(T^{m_k}x_0,T^mx_0)\le g(x_k) \circ...\circ g(x_0)(t_0),\  
m\ge m_k;
$$
wherefrom, taking (e09) into account,
$(T^nx_0; n\in N)$ is an increasing Cauchy sequence. 
Finally, given $y_0 \in Y$ with $x_0\le y_0$, put
$y_1= T^{n_0}y_0$ and, inductively, $y_{i+1}=T^{n_i}y_i=T^{m_i}y_0$, $i\ge 1$.
Again by (\ref{501}), 
$$
\mbox{
$d(x_0,T^mx_0),d(x_0,T^my_0)\le t_0$, $m\in N$, for some   $t_0> 0$.
}
$$
This fact, combined with (e07), leads us, by the same procedure as before, 
at 
$$
d(x_{k+1},T^mx_{k+1}),d(x_{k+1},T^my_{k+1})\le g(x_k) \circ ... \circ g(x_0)(t_0),\ m,k\in N;
$$
or equivalently (for each $k\in N$)
$$
d(T^{m_k}x_0,T^mx_0),d(T^{m_k}x_0,T^my_0)\le  g(x_k) \circ ... \circ g(x_0)(t_0),\ 
m\ge m_k;
$$
proving the desired conclusion and completing the argument.    
\end{proof}

{\bf (C)}
Let $X$, $d$ and  $\le$  be endowed with their previous meaning.
Given the sequence $(x_n; n\in N)$ in $X$ and the point $x\in X$,
define $x_n\uto x$ as: $(x_n; n\in N)$ is increasing and convergent to  $x$.   
Term the triplet $(X,d;\le)$, {\it quasi-order complete}, provided
each increasing Cauchy sequence converges. 
Note that any complete metric space
is quasi-order complete; but the converse is not in general valid. 
Further, given the self-map $T$ of $X$, call it {\it continuous at the left}
when $x_n \uto x$ and $x_n \le x$, $n \in N$, imply $Tx_n \to Tx$. 
Also, the ambient quasi-ordering  $\le$  will be said to be {\it self-closed} 
when $x\le y_n$, $n\in N$ and $y_n \uto y$ imply $x\le y$; 
note that any semi-closed quasi-ordering in 
Nachbin's sense \cite[Appendix]{nachbin-1965} is necessarily self-closed.
The first main result of the present note is

% Theorem 6.
\begin{theorem} \label{t6}
Let the conditions of Proposition \ref{p2} be fulfilled;
and (in addition)  
$(X,d;\le)$ is quasi-order complete,
$\le$  is self-closed, and
$T$ is continuous at the left.
Then, the following conclusions will be valid

{\bf iii)}
$Z:=\{x\in X; x=Tx\}$ is not empty

{\bf iv)}
for every $x\in Y$, $(T^nx; n\in N)$ converges to an element of $Z$

{\bf v)}
if $x,y\in Y$ are comparable, $(T^nx; n\in N)$ and $(T^ny; n\in N)$ 
have the same limit (in $Z$).
\end{theorem}

\begin{proof}
By Proposition \ref{p2} and the quasi-order completeness of $(X,d;\le)$,
it follows that, for the arbitrary fixed $x\in Y$, 
$T^nx\uto z$ for some $x\in X$. 
As $\le$ is self-closed, $T^nx \le z$, $n\in N$;
so that, combining with the left continuity of $T$ one gets 
$T^nx \uto Tz$; hence $z=Tz$.
The proof is thereby complete.
\end{proof}

Now, it is natural to ask of what happens when
$T$ is no longer continuous at the left. Some conventions are in order.
Call $\le$,  {\it anti self-closed} 
when $y_n \le x$, $n\in N$, and $y_n\uto y$ imply $y\le x$; 
observe at this moment that a sufficient condition for $\le$ to be anti self-closed 
is that  $\ge$ (its dual) be semi-closed. 
Further, call $\le$, {\it interval closed} 
when it is both self-closed and anti self-closed.
Our second main result is

% Theorem 7.
\begin{theorem}  \label{t7}
Let the conditions of Proposition \ref{p2} be fulfilled;
and (in addition)  
$(X,d;\le)$ is quasi-order complete
and $\le$ is an interval closed ordering.
Then, conclusions {\bf iii)}- {\bf v)} of Theorem \ref{t6} continue to hold;
and, moreover,

{\bf vi)}
for each $x\in Y$ the element $z=\lim_nT^nx$ in $Z$  has the properties 
(a)\ $x\le z$, (b)  $z\le y\in Y$ implies $z=y$.
\end{theorem}

\begin{proof}
Let $x\in Y$ be arbitrary fixed. 
By  Proposition \ref{p2}, $T^nx\uto z$, for some $z\in X$. 
hence (as $\le$ is self-closed), $x\le T^nx\le z$, $n\in N$.
It immediately follows that $T^nx\le Tz$, $n\in N$; 
so (by the anti self-closedness of $\le$), $z\in Y$.
Now, $x\le z\in Y$ gives, again by Proposition \ref{p2}, 
$T^nz \uto z$ (hence $Tz\le T^nz\le z$, $n\in N_0$)
and therefore (as $\le$ is ordering) $z\in Z$. The remaining part is evident.     
\end{proof}

It remains now to discuss the alternative: 
\bit
\item[]
($T$ is not continuous at the left) and 
($\le$  is not an interval closed ordering). 
\eit
To this end, assume that, for any $x\in Y$, 
the function $f(x)\in \calf_i(R_+^{2n(x)+1})$ given by (e07) fulfills
\bit
\item[(e10)]
for each $(\al_1,..,\al_{n(x)})\in R_+^{n(x)}$	with $\al_{n(x)}> 0$	
there exists $\be> 0$ with $\be+f(x)(\al_1,...,\al_{n(x)};\be,...,\be)< \al_{n(x)}$
\item[(e11)]
for each $(\al_1,...,\al_{n(x)})\in R_+^{n(x)}$ with $\al_1> 0$, $\al_{n(x)}=0$,
we have \\
$f(x)(\al_1,...,\al_{n(x)};\al_1,...,\al_{n(x)},\al_1)< \al_1$. 
\eit

Now, as a completion of the above results, we have

% Theorem 8.
\begin{theorem} \label{t8}
Let the conditions of Proposition \ref{p2} be fulfilled;
and (in addition)  
$(X,d;\le)$ is quasi-order complete,
(e10)+(e11) hold,
and $\le$ is an interval closed quasi-ordering.
Then, conclusions {\bf iii)}--{\bf vi)} of Theorem \ref{t7} 
still remain valid.
\end{theorem}

\begin{proof}
Let $x\in Y$ be arbitrary fixed. By the above reasoning,
$T^nx \uto z$, for some $z\in Y$;
with, in addition (cf. Proposition \ref{p2}):
[$x\le T^nx\le z$, $n\in N$] and  [$T^nz \uto z$]. 
Assume that $z\ne T^{n(z)}z$; and let 
$\be> 0$ be the number attached (via (e10)) to 
$\al_1:=d(z,Tz)$,..., $\al_{n(x)}:=d(z,T^{n(z)}z)$.
By the convergence property above, there exists $k(\be)\in N$ such that 
$d(z,T^kz)\le \be$, $\forall k\ge k(\be)$; and this gives
for all ranks $m\ge k(\be)+n(z)$,
$$  \barr{l}
d(z,T^{n(z)}z)\le  d(z,T^mz)+ d(T^{n(z)}z,T^mz)\le d(z,T^mz)+ \\
f(z)(d(z,Tz),...,d(z,T^{n(z)}z);d(z,T^{m-n(z)}z,...,d(z,T^mz))\le \\
\be+f(z)(d(z,Tz),...,d(z,T^{n(z)}z);\be,...,\be)< d(z,T^{n(z)}z); 
\earr
$$
contradiction; hence $z=T^{n(z)}z$. Moreover,
$$ \barr{l}
d(z,Tz)=d(T^{n(z)}z,T^{n(z)}Tz) \le \\
f(z)(d(z,Tz),...,d(z,T^{n(z)}z);d(z,Tz),...,d(z,T^{n(z)}z),d(z,T^{n(z)}Tz))= \\
f(z)(d(z,Tz),...,d(z,T^{n(z)-1}z),0;d(z,Tz),...,d(z,T^{n(z)-1}z),0,d(z,Tz));
\earr
$$
wherefrom, if $z\ne Tz$, (e11) will be contradicted.
Hence the conclusion.
\end{proof}

Some remarks are in order.
Theorem \ref{t6} may he viewed as a quasi-order extension 
of Sehgal's result we just quoted  
(cf. also Dugundji and Granas \cite[Ch 1, Sect 3]{dugundji-granas-1982})
while Theorem \ref{t7} is a quasi-order "functional" version of
Matkowski's contribution (Theorem \ref{t5}). 
At the same time, Theorem \ref{t8} - although formulated as
a fixed point result - may be deemed in fact as a maximality principle 
in $(Y,\le)$; so, it is comparable under this perspective with a related 
author's one \cite{turinici-1981} 
obtained by means of a "compactness" procedure like in
Krasnoselskii and Sobolev \cite{krasnoselskii-sobolev-1973}.
\bk

{\bf (D)}
{\it Note added in 2011}
\sk

From these developments,
the following statement is deductible.
Let  the quasi-ordered metric space $(X,d,\le)$ 
the self-map $T$ of $X$ be taken as in (e05)+(e06).
In addition, the specific condition will be accepted:
\bit
\item[(e12)]
there exists $f\in \calf_i(R_+)$ such that: 
for each $x$ in $Y$ there exists $n(x)\in N_0$ \\
with
$d(T^{n(x)}x,T^{n(x)}y)\le f(d(x,y))$,
for all $y\in Y$ with $x\le  y$.
\eit
Note that, in such a case, the iterative normality of 
$((n(x);f); x\in Y)$ is characterized by (e02)+(e03),
with $f$ in place of $g$;
and referred to as: $f$ is {\it normal} (see above).
From Theorem \ref{t6} we then get, formally

% Theorem 9.
\begin{theorem} \label{t9}
In addition to (e05)+(e06), assume that the function 
$f$ (appearing in (e12)) is normal, 
$(X,d;\le)$ is quasi-order complete,
$\le$  is self-closed, and
$T$ is continuous at the left.
Then, conclusions {\bf iii)}--{\bf v)} of Theorem \ref{t6}
are retainable.
\end{theorem}

In particular, any linear comparison function $f$
(in the sense: $f(t)=\al t$, $t\in R_+$, for $0< \al< 1$)
is normal. 
Then, Theorem \ref{t9} includes the essential conclusions
of the Ran-Reurings result (Theorem \ref{t1}).
[In fact, under appropriate conditions, it may give us all conclusions 
in the quoted statement; we do not give details].
Note that Theorem \ref{t9} is not yet covered by the
existing fixed point statements in the 
realm of quasi-ordered metric spaces.
Further aspects will be delineated elsewhere.

%  References


\begin{thebibliography}{99}


% 1
\bibitem{agarwal-el-gebeily-o-regan-2008}
{R. P. Agarwal}, {M. A. El-Gebeily} and {D. O'Regan},
\it Generalized contractions in partially ordered metric spaces, 
\rm Appl. Anal., 87 (2008), 109-116. 


% 2
\bibitem{banach-1922}
{S. Banach}, 
\it Sur les op\'{e}rations dans les ensembles abstraits 
et leur application aux \'{e}quations int\'{e}grales.
\rm Fund. Math., 3 (1922), 133-181.


% 3
\bibitem{ciric-1972}
{L. B. Ciric},     
\it Fixed point theorems for mappings with a generalized contractive iterate at a point, 
\rm Publ. Inst. Math., 13 (27) (1972), 11-16. 

% 4
\bibitem{ciric-mihet-saadati-2009}
{L. B. Ciric}, {D. Mihet} and {R. Saadati},
\it Monotone generalized contractions in partially ordered probabilistic metric spaces,
\rm Topology and its Appl., 156 (2009), 2838-2844.

% 5
\bibitem{dugundji-granas-1982}
{J. Dugundji} and {A. Granas},    
\it Fixed Point Theory, 
\rm vol. I, Warszawa, 1982. 


% 6
\bibitem{frechet-1906}
{M. Fr\'{e}chet},
\it Sur quelques points du calcul fonctionnel,
\rm Rend. Circ. Mat. Palermo, 22 (1906), 1-72.


% 7
\bibitem{guseman-1970}
{L. F. Guseman Jr.},   
\it Fixed point theorems for mappings with a contractive iterate at a point, 
\rm Proc. Amer. Math. Soc., 26 (1970), 615-618. 


% 8
\bibitem{gwozdz-lukawska-jachymski-2009}
{G. Gwozdz-Lukawska} and {J. Jachymski},
\it IFS on a metric space with a graph structure and extensions of the Kelisky-Rivlin theorem, 
\rm J. Math. Anal. Appl., 356 (2009), 453-463. 


% 9
\bibitem{iseki-1974}
{K. Iseki},     
\it A generalization of Sehgal-Khazanchi's fixed point theorem, 
\rm Math. Sem. Notes Kobe Univ., 2 (1974), 89-95.


% 10
\bibitem{jachymski-1993}
{J. Jachymski}, 
\it A generalization of the theorem by Rhoades and Watson for contractive type mappings, 
\rm Math. Japon. 38 (1993), 1095-1102. 


% 11
\bibitem{jachymski-1994}
{J. Jachymski},
\it Common fixed point theorems for some families of mappings,
\rm  Indian J. Pure Appl. Math., 25 (1994), 925-937.


% 12
\bibitem{kasahara-1976}
{S. Kasahara},
\it On some generalizations of the Banach contraction theorem,
\rm Publ. Res. Inst. Math. Sci. Kyoto Univ., 12 (1976), 427-437.


% 13
\bibitem{khazanchi-1974}
{L. Khazanchi},     
\it Results on fixed points in complete metric space, 
\rm Math. Japon., 19 (1974), 283-289.


% 14
\bibitem{krasnoselskii-sobolev-1973}
{M. A. Krasnoselskii} and {A. V. Sobo1ev}, 
\it O nepodvizhnych tochkach razryvnych operatorov, 
\rm  Sibirsk. Mat. Zh., 14 (1973) 674-677.


% 15
\bibitem{maia-1968}
{M. G. Maia},
\it Un'osservazione sulle contrazioni metriche,
\rm Rend. Sem. Mat. Univ. Padova, 40 (1968), 139-143.



% 16
\bibitem{matkowski-1977}
{J. Matkowski},     
\it Fixed point theorems for mappings with a contractive iterate at a point, 
\rm Proc. Amer. Math. Soc., 62 (1977), 344-348.


% 17
\bibitem{nachbin-1965}
{L. Nachbin}, 
\it Topology and Order, 
\rm Van Nostrand, Princeton, N.J., 1965.


% 18
\bibitem{nieto-rodriguez-lopez-2005}
{J.J. Nieto} and  {R. Rodriguez-Lopez}, 
\it Contractive mapping theorems in partially ordered sets and applications 
to ordinary differential equations,
\rm Order, 22 (2005), 223-239.


% 19
\bibitem{o-regan-petrusel-2008}
{D. O'Regan} and  {A. Petru\c{s}el}, 
\it Fixed point theorems for generalized contractions in ordered metric spaces, 
\rm J. Math. Anal. Appl., 341 (2008), 1241-1252.


% 20
\bibitem{petrusel-rus-2006}
{A. Petru\c{s}el} and  {I. A. Rus}, 
\it Fixed point theorems in ordered $L$-spaces, 
\rm Proc. Amer. Math. Soc., 134 (2006), 411-418.


% 21
\bibitem{ran-reurings-2004}
{A. C. M. Ran} and  {M. C. Reurings}, 
\it A fixed point theorem in partially ordered sets and some applications to matrix equations, 
\rm Proc. Amer. Math. Soc., 132 (2004), 1435-1443.


% 22
\bibitem{rhoades-1977}
{B. E. Rhoades}, 
\it A comparison of various definitions of contractive mappings, 
\rm Trans. Amer. Math. Soc., 226 (1977), 257-290.


% 23
\bibitem{rus-2001}
{I. A. Rus}, 
\it Generalized Contractions and Applications, 
\rm Cluj University Press, Cluj-Napoca, 2001.


% 24
\bibitem{sehgal-1969}
{V. M. Sehgal},      
\it A fixed point theorem for mappings with a contractive iterate, 
\rm Proc. Amer. Math. Soc., 23 (1969), 631-634.


% 25
\bibitem{singh-1979}
{K. L. Singh}, 
\it Fixed point theorems for contractive-type mappings, 
\rm J. Math. Anal. Appl., 72 (1979), 283-290.


% 26
\bibitem{turinici-1981}
{M. Turinici},     
\it A class of operator equations on ordered metric spaces, 
\rm Bull. Malaysian Math. Soc., (2), 4 (1981), 67-72.


% 27
\bibitem{turinici-1986}
{M. Turinici},
\it Fixed points for monotone iteratively local contractions,
\rm Dem. Math., 19 (1986), 171-180.



\end{thebibliography}
\end{document}